\documentclass[english]{article}
\pdfoutput=1

\usepackage[T1]{fontenc}
\usepackage[latin9]{inputenc}
\usepackage{verbatim}
\usepackage{float}
\usepackage{amsmath}
\usepackage{ifpdf}
\usepackage{times}
\usepackage{pifont}

\usepackage[margin=10pt,font=small,labelfont=bf,labelsep=endash]{caption}

\ifpdf
  \usepackage[pdftex]{graphicx}
  \usepackage[pdftex,dvipsnames,svgnames]{xcolor}
\else
  \usepackage[dvips,dvipsnames,svgnames]{xcolor}
  \usepackage[dvips]{graphicx}
\fi

\usepackage{wrapfig}
\usepackage{subfig}
\usepackage{setspace}
\usepackage[authoryear]{natbib}

\usepackage{babel}
\makeatletter

\@ifundefined{c@subfigure}{\newsubfloat{figure}}{}
\def\subfigure{\subfloat}

\begin{document}

\title{Understanding Terrorist Organizations with a Dynamic Model}

\author{Alexander Gutfraind}


\maketitle
\begin{abstract}
Terrorist organizations change over time because of processes such
as recruitment and training as well as counter-terrorism (CT) measures,
but the effects of these processes are typically studied qualitatively
and in separation from each other. Seeking a more quantitative and
integrated understanding, we constructed a simple dynamic model where
equations describe how these processes change an organization's membership.
Analysis of the model yields a number of intuitive as well as novel
findings. Most importantly it becomes possible to predict whether
counter-terrorism measures would be sufficient to defeat the organization.
Furthermore, we can prove in general that an organization would collapse
if its strength and its pool of foot soldiers decline simultaneously.
In contrast, a simultaneous decline in its strength and its pool of
leaders is often insufficient and short-termed. These results and
other like them demonstrate the great potential of dynamic models
for informing terrorism scholarship and counter-terrorism policy making.
\end{abstract}

\section{Introduction}

Our goal is to study terrorist organizations using a dynamic model.
Generally speaking, in such a model a phenomenon is represented as
a set of equations which describe it in simplified terms. The equations
represent how the phenomenon changes in time or space, and cast our
empirically-based knowledge in precise mathematical language. Once
the model is constructed, it can be studied using powerful mathematical
techniques to yield predictions, observations and insights that are
difficult or impossible to collect empirically \citep{Aris95,Ellner06}.
 For example, a dynamic model could be constructed for the various
militant groups operating in Iraq and then used to predict their strength
a year in the future. Moreover, given the model, it would be possible
to evaluate the efficacy of various counter-insurgency policies.

Mathematical models can help fill a large methodological void in terrorism
research: the lack of model systems. Whereas biologists studying pathogens
can do experiments \emph{in vitro}, there are no such model systems
in terrorism research, except for mathematical models. In this sense,
the method developed below provides an \emph{in vitro} form of terrorism,
which can be investigated in ways not possible in its \emph{in vivo}
kind. Like all model systems, mathematical models are imperfect because
they rely on large simplifications of the underlying political phenomena,
and one can rightfully ask whether their predictions would be sufficiently
accurate. Fortunately, complex phenomena in fields like biology have
been studied very successfully with this mathematical technique \citep{Ellner06}.
Therefore, even phenomena as complex as found in terrorism research
may, in some cases, be productively studied using mathematical models
and indeed, previous models have brought considerable insights%
\footnote{E.g. dynamic models: \citet{Allanach04,Chamberlain07,Farley07,Feichtinger01,Johnson05,Stauffer06,Udwadia06},
rational choice models: \citet{Anderton05,Sandler83,Sandler03,Wintrobe06},
agent-based models: \citet{MacKerrow03,Tsvetovat07}.%
}.  

In the rest of the paper we describe a simple model of a terrorist
organization. The model is new in its focus, methodology and audience:
We focus on a single terrorist organization and model its processes
of recruitment, its internal dynamics as well as the impact of counter-terrorism
measures on it. As to methodology, with a few exceptions \citep{Chamberlain07,Feichtinger01,Stauffer06,Udwadia06}
and perhaps a few others the powerful mathematical technique of differential equations has
not been applied in terrorism research. Finally, the paper
is specifically written to an audience of non-mathematicians: the
main body of the paper uses non-technical language to explain the
terminology and to describe the equations and assumptions used in
the model, while the technical analysis is exposed in the appendix.

The model described below was built following two design principles.
First, it was desired to have a model of broad applicability across
organizations and conflicts. Indeed, the model is so general that
it can be applied to insurgencies or even to some non-terrorist organizations.
As we shall see, despite this generality it makes non-trivial observations and
more importantly it specifies sufficient conditions for victory over
the organization (see subsection \ref{sub:Conditions-for-Victory}).
Second, it was desired to build a simple model so as to facilitate
interpretation, analysis and further development. It was hoped that
the model would establish a methodological prototype that could be
easily extended and modified to fit specific cases.

The organization of the paper is as follows. Section 2 describes the
model - its variables, parameters and relations between them. Section
3 graphically illustrates the model's predictions about terrorist
organizations. Sections 4 and 5 discuss the insights gleaned from
the model, and the implications to counter-terrorism policies. The
conclusions are in Section 6. Finally, all of the technical arguments
are gathered in the appendix.

\section{\label{sec:A-Mathematical-Model}A Mathematical Model}

There are many ways of describing a terrorist organization, such as
its ideology or political platform, its operational patterns, or its
methods of recruitment. Here we consider it from the {}``human resources''
point of view. Namely, we are interested in examining how the numbers
of {}``leaders'' and {}``foot soldiers'' in the organization change
with time. The former includes experienced managers, weapon experts,
financiers and even politicians and religious leaders who help the
organization, while the latter are the more numerous rank-and-file.
These two quantities arguably give the most important information
about the strength of the organization. The precise characteristics
of the two groups and their relative sizes would depend on the organization
under consideration. Nevertheless, this distinction remains relevant
even in the very decentralized organizations like the post-Afghanistan
al-Qaeda movement, because we can identify the ``leaders'' as
the experienced terrorists, as compared to the new recruits 
(see discussions in \citealp{Hoffman03,Sageman04}).
The division between those two groups is also important in practice
because decision makers often need to choose which of the groups
to target \citep[Ch.5]{Wolf89,Ganor05}: while the leaders
represent more valuable targets, they are also harder to reach. Later
on in section 5 we actually compare the policy alternatives.

Therefore, let us represent a terrorist organization as two time-varying
quantities, $L$ and $F$, corresponding to the number of leaders
and foot soldiers, respectively. Also, $L$ and $F$ determine the
overall {}``strength'' $S$ of the organization. Because
leaders possess valuable skills and experience, they contribute more
to the strength than an equivalent number of foot soldiers. Hence,
strength $S$ is taken to be a weighted sum of the two variables,
with more weight ($m>1$) given to leaders:\[
S=mL+F\]

We now identify a set of processes that are fundamental in 
changing the numbers of leaders and foot soldiers. These
processes constitute the mathematical model. While some of them
are self-evident, others could benefit from quantitative comparison
with data. The latter task is non-trivial given the scarcity of time-series
data on membership in terrorist organizations and hence we leave it
out for future work.

The histories of al-Qaeda and other terrorist organizations \citep[e.g.][]{Laqueur01,Harmon00,Hoffman06}
suggest that the pool of terrorist leaders and experts grows primarily
when foot soldiers acquire battle experience or receive training (internally,
or in terrorist-supporting states, \citealp{Jongman05}). Consequently,
the pool of leaders ($L$) is provisioned with new leaders at a rate
proportional to the number of foot soldiers ($F$). We call this process
{}``promotion'' and label the parameter of proportionality $p$.
This growth is opposed by internal personnel loss due to demotivation,
fatigue, desertion as well as in-fighting and splintering \citep[Ch.6]{Horgan05}.
This phenomenon is modeled as a loss of a fraction ($d$) of the pool
of leaders per unit time. An additional important influence on the
organization are the counter-terrorism (CT) measures targeted specifically
at the leadership, including arrests, assassinations as well as efforts
to disrupt communications and to force the leaders into long-term
inactivity. Such measures could be modeled as the removal of a certain
number ($b$) of people per unit time from the pool of leaders. CT
is modeled as a constant rate of removal rather than as a quantity
that depends on the size of the organization because the goal is to
see how a fixed resource allocation towards CT would impact the organization.
Presumably, the human resources and funds available to fight the given
terrorist organization lead, on average, to the capture or elimination
of a fixed number of operatives. In sum, we assume that on average,
at every interval of time the pool of leaders is nourished through promotion,
and drained because of internal losses and CT (see appendix, equation
(\ref{eq:L}).)

The dynamics of the pool of foot soldiers ($F$) are somewhat similar
to the dynamics of leaders. Like in the case of leaders, some internal
losses are expected. This is modeled as the removal of a fraction
($d$) of the pool of operatives per unit time where for simplicity
the rate $d$ is the same as the rate for leaders (the complex case is discussed 
in subsection \ref{sub:Encouraging-desertion}.) Much like in the
case of leaders above, counter-terrorism measures are assumed to remove
a fixed number ($k$) of foot soldiers per unit time. Finally and
most importantly, we need to consider how and why new recruits join
a terrorist organization. Arguably, in many organizations growth in
the ranks is in proportion to the strength of the organization, for
multiple reasons: Strength determines the ability to carry out successful
operations, which increase interest in the organization and its mission.
Moreover, strength gives the organization the manpower to publicize
its attacks, as well as to man and fund recruitment activities. By
assuming that recruitment is proportional to strength, we capture
the often-seen cycle where successful attacks lead to greater recruitment,
which leads to greater strength and more attacks. Overall, the pool
of foot soldiers is nourished through recruitment, and drained because
of internal losses and CT (see appendix, equation (\ref{eq:F}))%
\footnote{A minor assumption in our model is that once a foot soldier is promoted
a new foot soldier is recruited as a replacement. It is shown in the
appendix that if in some organizations such recruitment is not automatic,
then the current model is still valid for these organizations as long
as $p<r$. In any case the drain due to promotion is marginal because
foot soldiers are far more numerous than leaders even in relatively
{}``top heavy'' organizations.%
}$^,$%
\footnote{This model is similar to structured population models in biology,
where the foot soldiers are the {}``larvae'' and the leaders are
the {}``adults''. However, an interesting difference is that whereas
larvae growth is a function of the adult population alone, in a terrorist
organization the pool of foot soldiers contributes to its own growth.%
}. 

The numerical values of all of the above parameters ($p,d,b,r,m,k$)
are dependent on the particular organization under consideration,
and likely vary somewhat with time%
\footnote{The simplest approach to estimating them would be to estimate the
number and leaders and foot soldiers at some point in time, and then
find the parameter values by doing least-squares fitting of the model
to the data on the terrorist attacks, where we consider the terrorist
attacks to be a proxy of strength. However, this approach has some
limitations.%
}. Fortunately, it is possible to draw many general conclusions from
the model without knowing the parameter values, and we shall do so
shortly. Finally, it should be noted that counter-terrorism need not
be restricted to the parameters $b,k$ (removal of leaders and foot
soldiers, respectively), and measures such as public advocacy, attacks
on terrorist bases, disruption of communication and others can weaken
the organization by reducing its capabilities as expressed through
the other parameters.

In the above description, we assumed that counter-terrorism measures
are parameters that can be changed without affecting recruitment.
This is a significant simplification because in practice it may be
difficult to respond to terrorist attacks without engendering a backlash
that actually promotes recruitment \citep[see e.g. ][]{Ganor08,Hanson07}.
Nevertheless, the advantages of this simplification outweigh the disadvantages:
Firstly, it is clear that any model that would consider such an effect
would be much more complicated than the current model and consequently
much harder to analyze or use. Secondly, the current model can be
easily extended to incorporate such an effect if desired. Thirdly,
the strength of this effect is difficult to describe in general because
it depends extensively on factors such as the the specific CT measures
being used, the terrorist actions and the political environment. Indeed,
\citet{Udwadia06}\emph{ }who incorporated this effect, constructed
their model based on observations of a specific context within the
current conflict in Iraq.

The model includes additional implicit assumptions. First,  it assumes
a state of stable gradual change, such that the effect of one terrorist
or counter-terrorism operation is smoothed. This should be acceptable
in all cases where the terrorist organization is not very small and
thus changes are not very stochastic. Second, the model assumes that
an organization's growth is constrained only by the available manpower,
and factors such as money or weapons do not impose an independent
constraint. Third, it is assumed that the growth in foot soldiers
is not constrained by the availability of potential recruits - and
it is probably true in the case of al-Qaeda because willing recruits
are plentiful (for the case of England, see \citealp{Manningham-Buller06}).
We discuss this point further in subsection \ref{sub:Stable-Equilibria}.

\section{Analysis of the Model}

Having written down the governing equations, the task of studying
a terrorist organization is reduced to the standard problem of studying
solutions to the equations. Because the equations indicate rates of
change in time, the solutions would be two functions, $L(t)$ and
$F(t)$, giving the number of leaders and foot soldiers, respectively,
at each time. Let us suppose that currently (time $0$) the organization
has a certain number of leaders and foot soldiers, $L_{0}$ and $F_{0}$,
respectively and is subject to certain CT measures, quantified by
$b$ and $k$. We want to see whether the CT measures are adequate
to defeat the organization. Mathematically, this corresponds to the
question of whether at some future time both $L$ and $F$ would reach
zero. Intuitively, we expect that the organization would be eliminated
if it is incapable of recovering from the losses inflicted on it by
CT. In turn, this would depend on its current capabilities as well
as the parameters $p,d,r,m$ which characterize the organization.

\begin{figure}[htb]
\begin{centering}
\subfigure[]{\includegraphics[width=0.5\columnwidth]{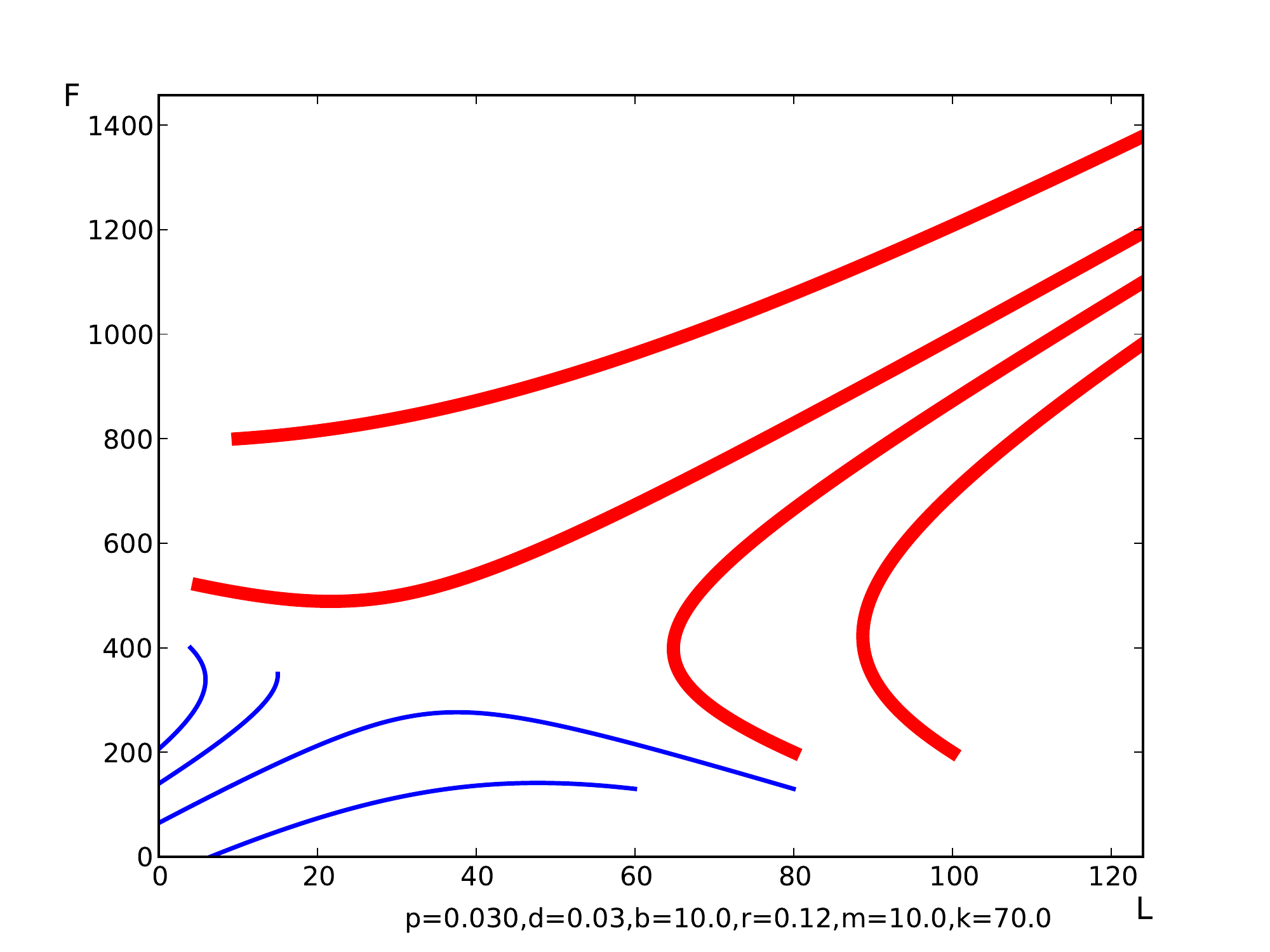}}\subfigure[]{\includegraphics[width=0.5\columnwidth]{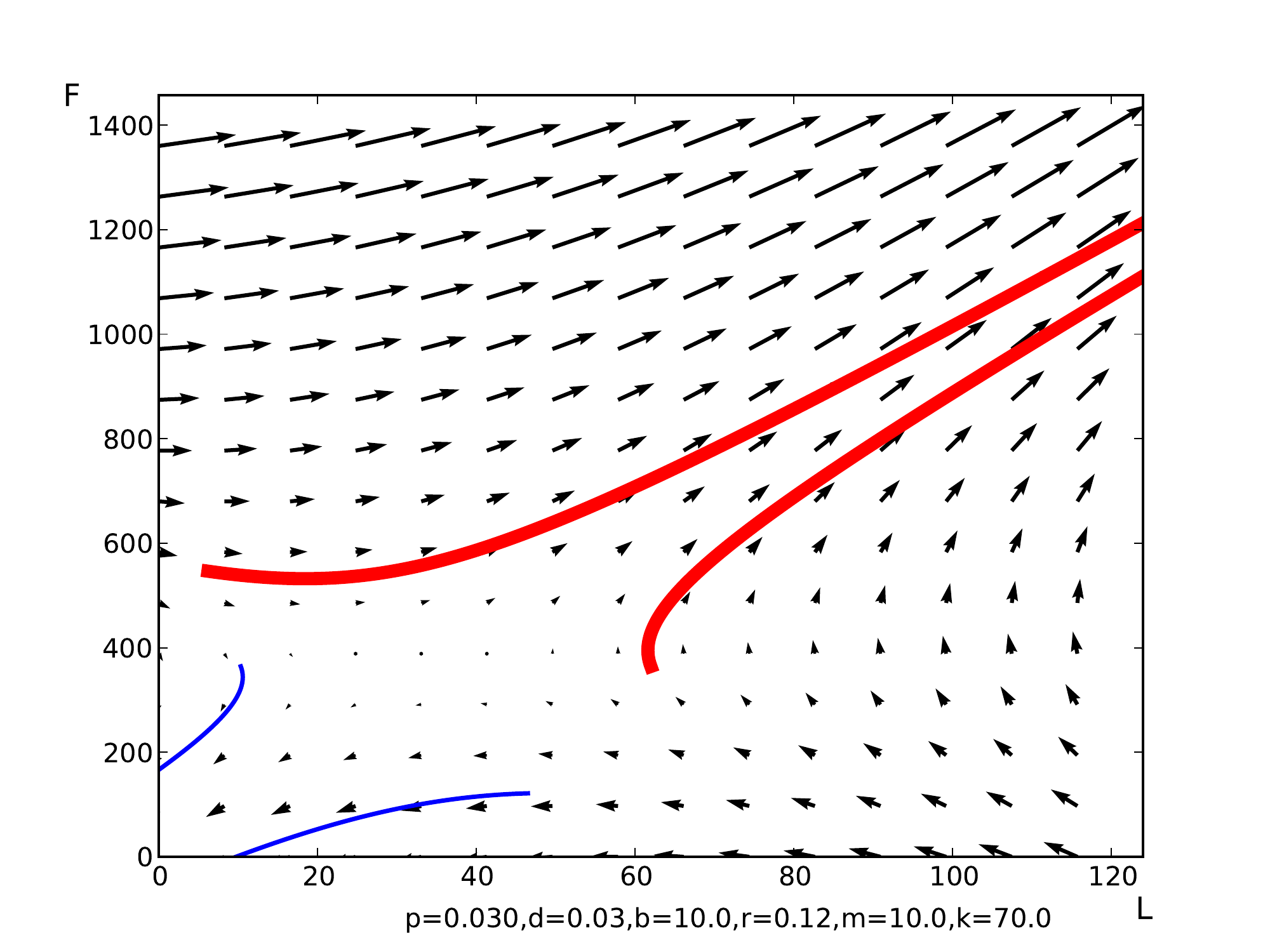}}
\caption{\label{fig:Phase-portrait}(a) Typical solution curves of the equations
coded by ultimate fate: thin blue for successfully neutralized organizations and
thick red for those remaining operational and growing. The parameters
were set to representative values, but as was said earlier, all realistic
organizations are qualitatively similar and resemble these. (b) {}``Vector
field'' of $L$ and $F$. At each value of $L,F$ the direction and
length of the arrow give the rate of change in $L$ and $F$. }
\par\end{centering}
\end{figure}

\begin{figure}[H]
\begin{centering}
\subfigure[]{\includegraphics[width=0.5\columnwidth]{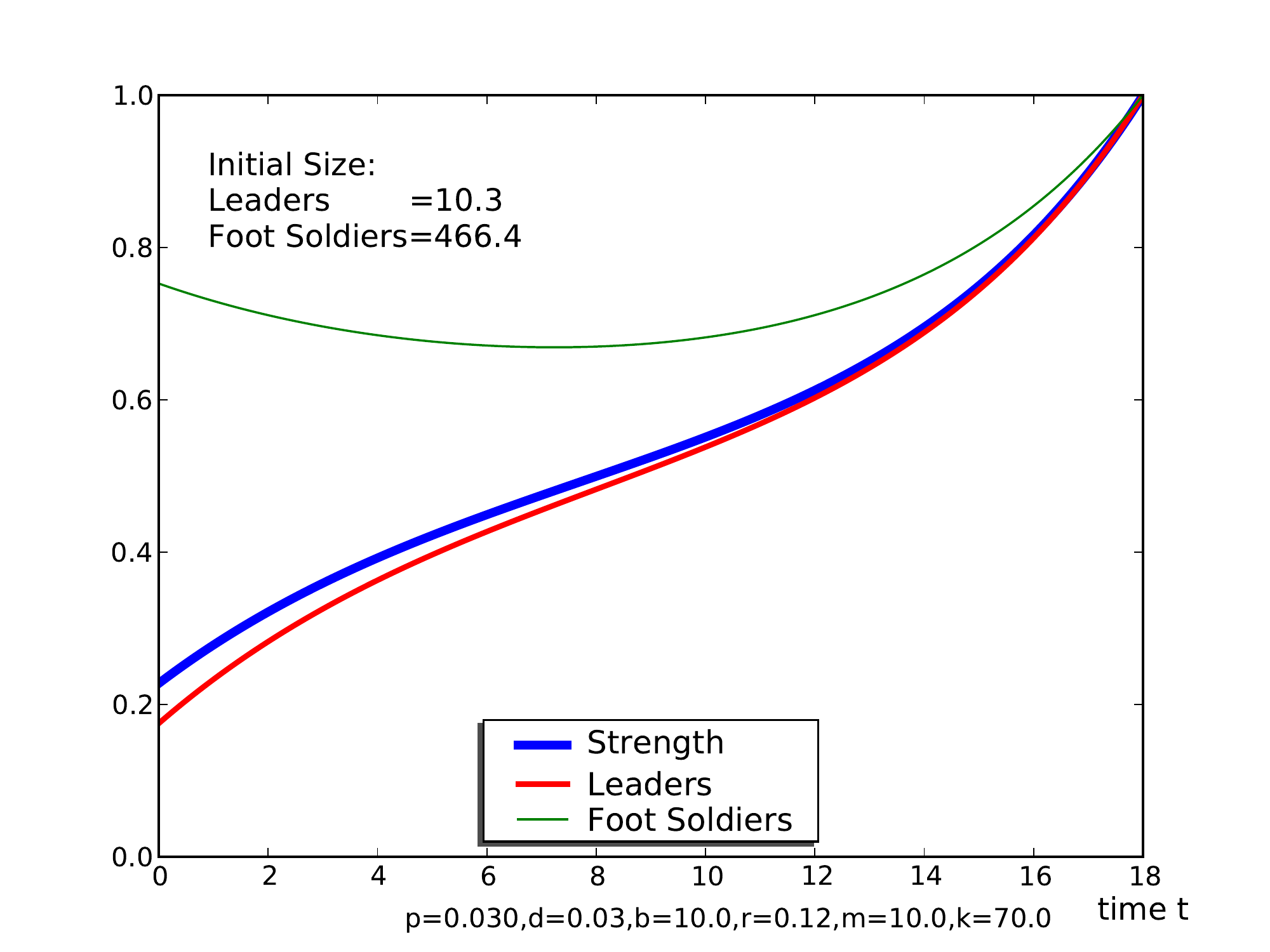}}\subfigure[]{\includegraphics[width=0.5\columnwidth]{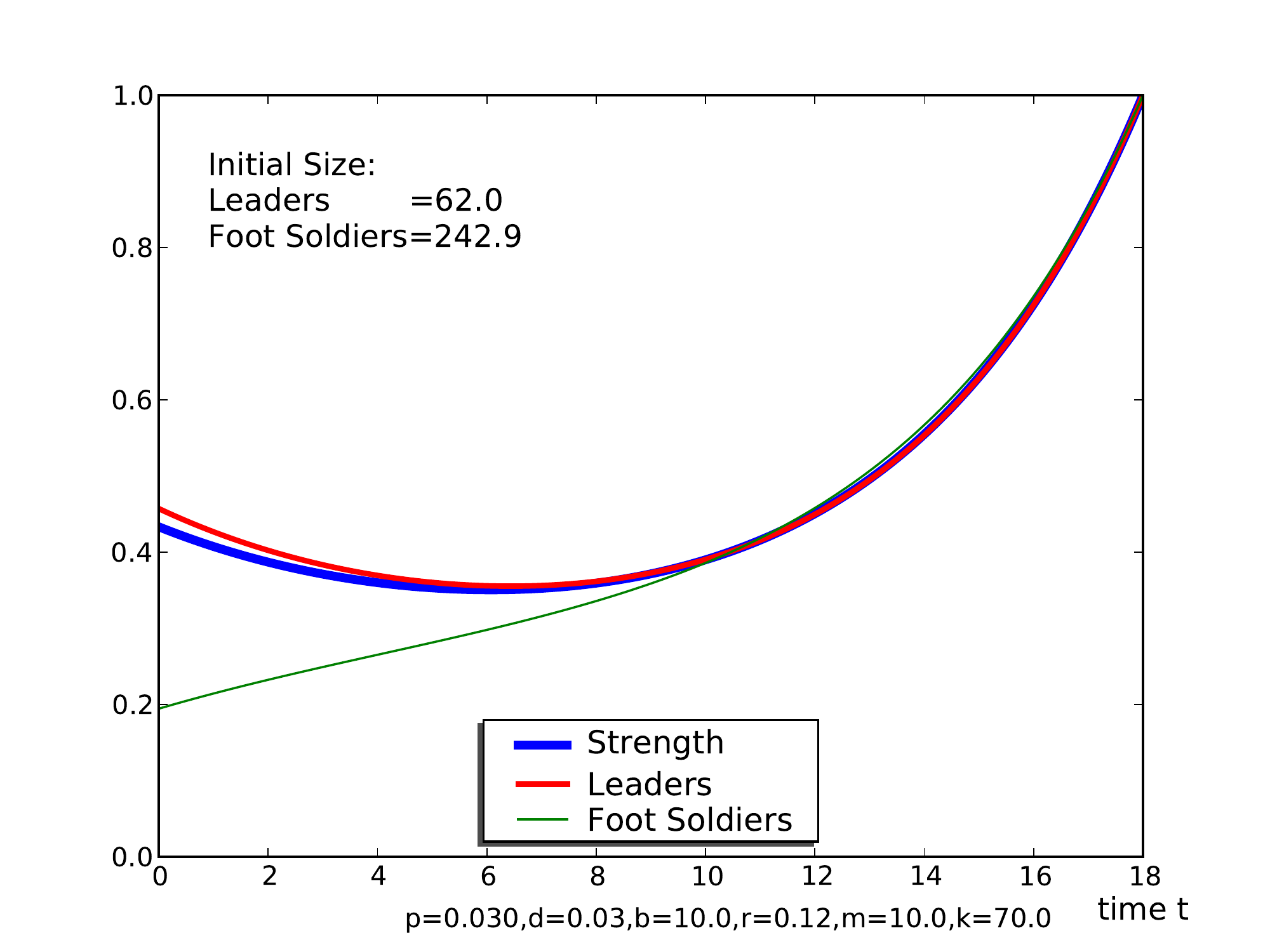}}
\par\end{centering}

\caption{\label{fig:Evolution1}Evolution of strength, leaders and foot soldiers
($S,L,F$, respectively) in two terrorist organizations as a function
of time. In (a), due to CT, $F$ falls initially but eventually the
organization recovers through promotion. In (b), $L$ and $S$ fall
initially but eventually the organization recovers through recruitment.
The vertical axis has been rescaled by dividing each quantity by the
maximum it attains during the time evolution. This makes it possible
to represent all three quantities on the same plot. The units of time
are unspecified since they do not affect the analysis. 
Of course, in a more complex model it would be desirable to consider 
periodic events like election cycles or generational changes.}

\end{figure}

Mathematical analysis of the model (see the appendix)
shows that most terrorist organizations%
\footnote{That is, those with realistically low rates of desertion: $d<\frac{1}{2}(r+r\sqrt{1+\frac{mp}{r}})$.
A higher rate of desertion $d$ always causes the organization to
collapse and is not as interesting from a policy perspective (see
subsection \ref{sub:Encouraging-desertion} for a discussion of desertion).%
} evolve in time like the organizations whose {}``orbits'' are displayed
in Fig.~\ref{fig:Phase-portrait}(a,b). In Fig.~\ref{fig:Phase-portrait}(a)
we plotted eight different organizations with different starting conditions.
Another perspective can be seen in Fig.~\ref{fig:Phase-portrait}(b)
which graphically illustrates the dynamical equations via arrows:
the direction of each arrow and its length indicate how an organization
at the tail of the arrow would be changing and it what rate. By picking
any starting location ($L_{0},F_{0}$) and connecting the arrows,
it is possible to visually predict the evolution into the future.
Another illustration is found in Fig.~\ref{fig:Evolution1}, which
shows how two example organizations change with time.

\begin{wrapfigure}{r}{0.5\textwidth}
\begin{centering}
\includegraphics[width=0.5\textwidth]{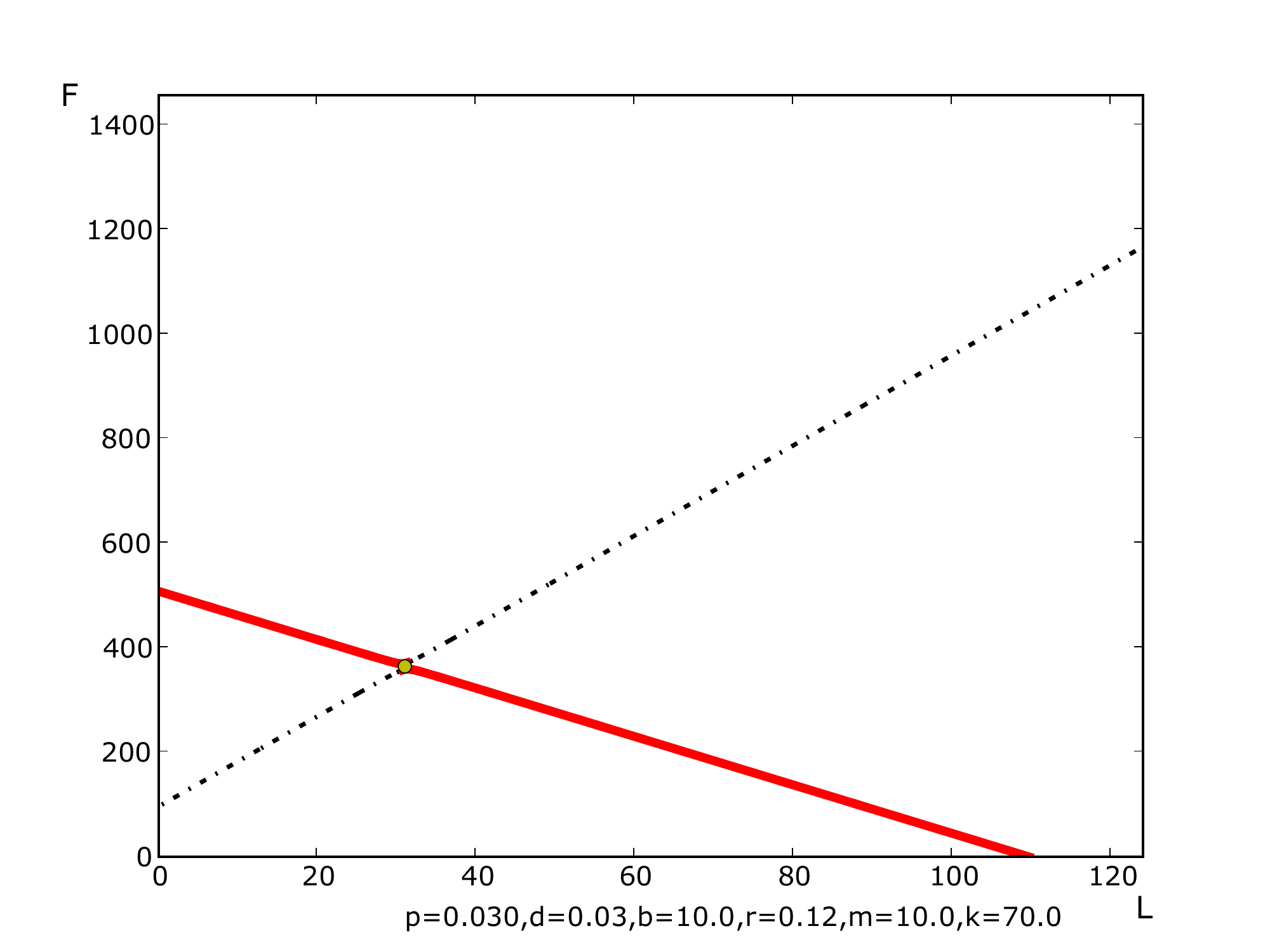}
\caption{\label{fig:Phase-Schematical}Plot of the sink (thick red) and trend
lines (thin dashed black). The two lines intersect at a {}``saddle
point''. 
}
\par\end{centering}
\end{wrapfigure}

In general, it is found that the dynamics of the organization is dependent
upon the position of the organization with respect to a threshold
line, which can be termed the {}``sink line'': an organization will
be neutralized if and only if its capabilities are below the sink
line. In other words, the current CT measures are sufficient if and
only if the organization lies below that threshold (thick red line
on Fig.~\ref{fig:Phase-Schematical}). The threshold is impassable:
an organization above it will grow, and one below it is sure to collapse.
This threshold is also very sharp: two organizations may lie close
to the line, but the one above it would grow, while the one below
it would shrink even if the differences in initial capabilities are
small. In addition to the sink line, the model also predicts that
all successful organizations would tend towards a particular trajectory.
This {}``trend line'' (a dashed black line on Fig.~\ref{fig:Phase-Schematical})
is discussed further in subsection \ref{sub:classification}.

\begin{center}
\begin{figure}[H]
\begin{centering}
\subfigure[]{\includegraphics[width=0.4\columnwidth]{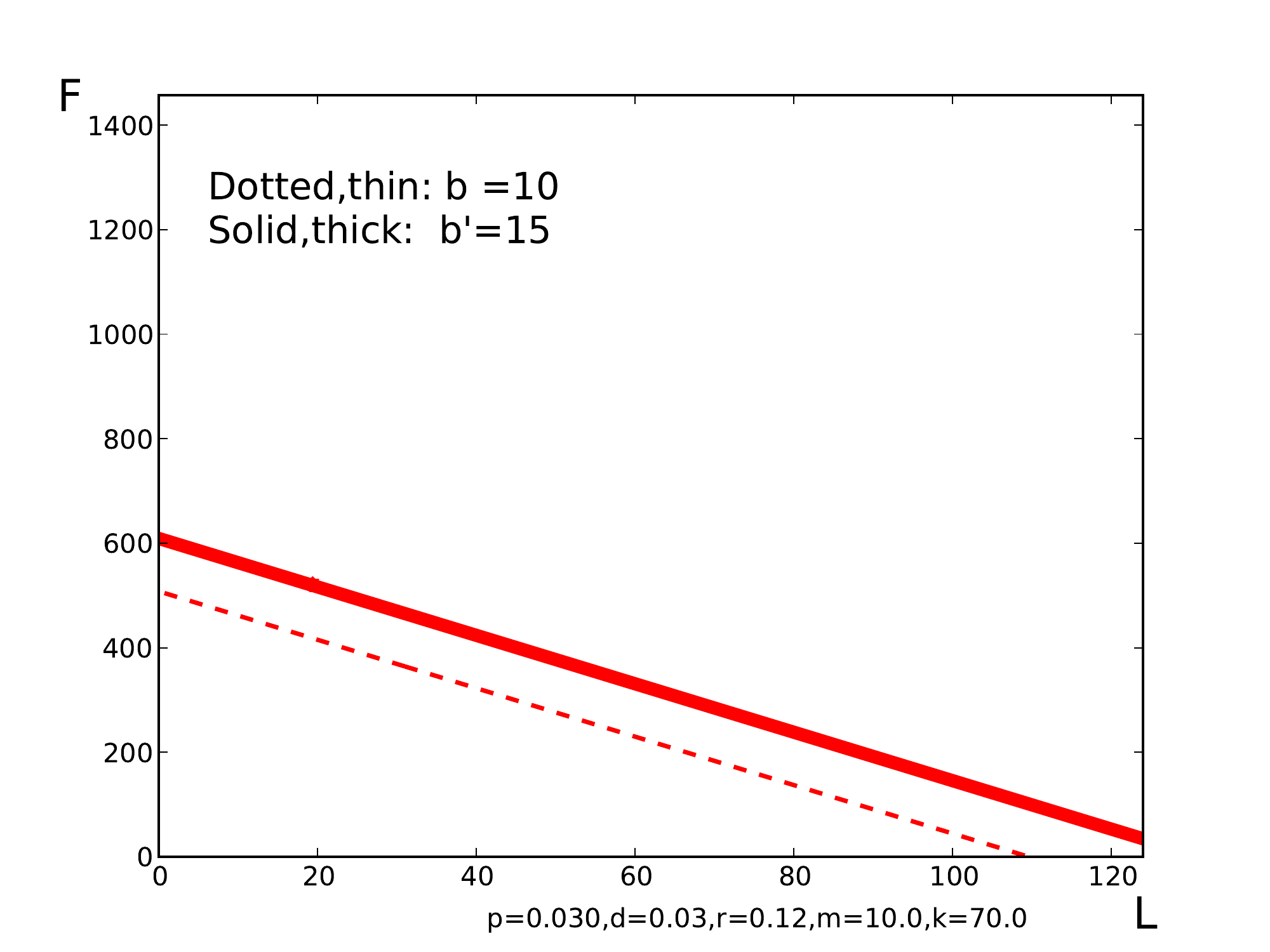}}\subfigure[]{\includegraphics[width=0.4\columnwidth]{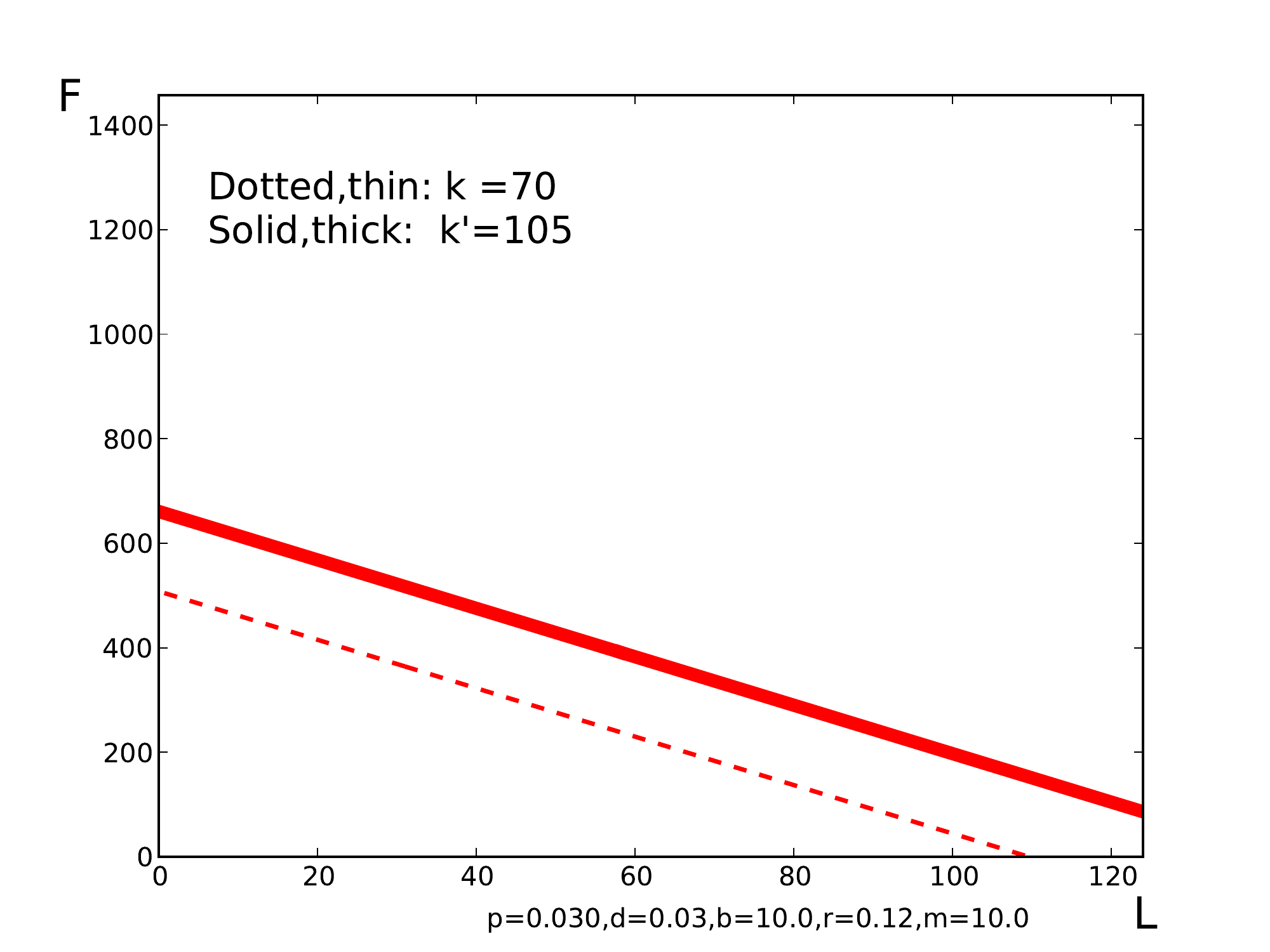}}
\par\end{centering}
\caption{\label{fig:The-effects-of-bk}The effects of the parameters $b$ and
$k$ on the dynamical system, (a) and (b) respectively, as seen through
the effect on the sink line. In each case, as the CT measures are
increased, the sink line moves up confining below it additional terrorist
organizations.}
\end{figure}
\end{center}

\begin{figure}[H]
\begin{centering}
\subfigure[]{\includegraphics[width=0.3\columnwidth]{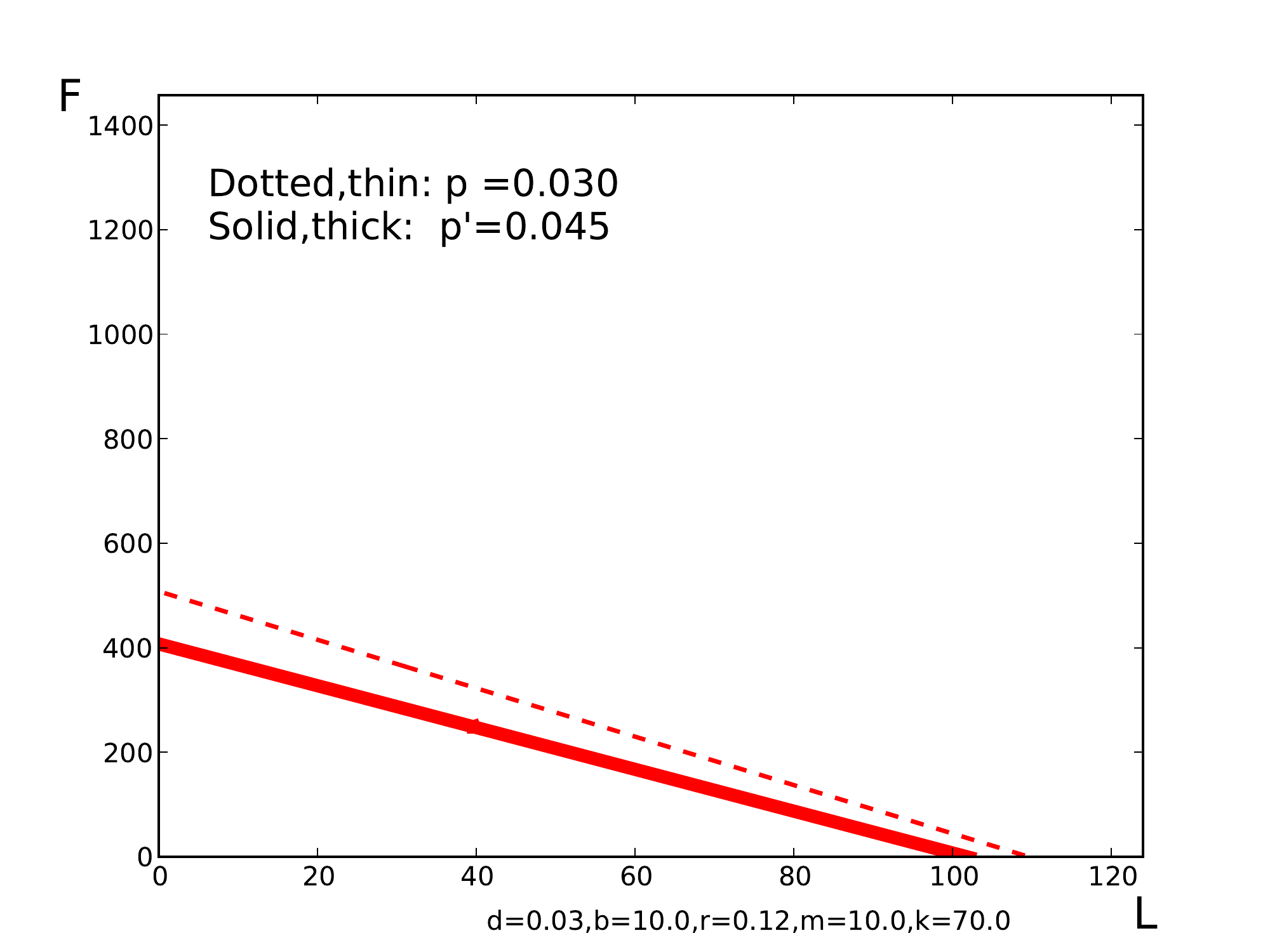}}\subfigure[]{\includegraphics[width=0.3\columnwidth]{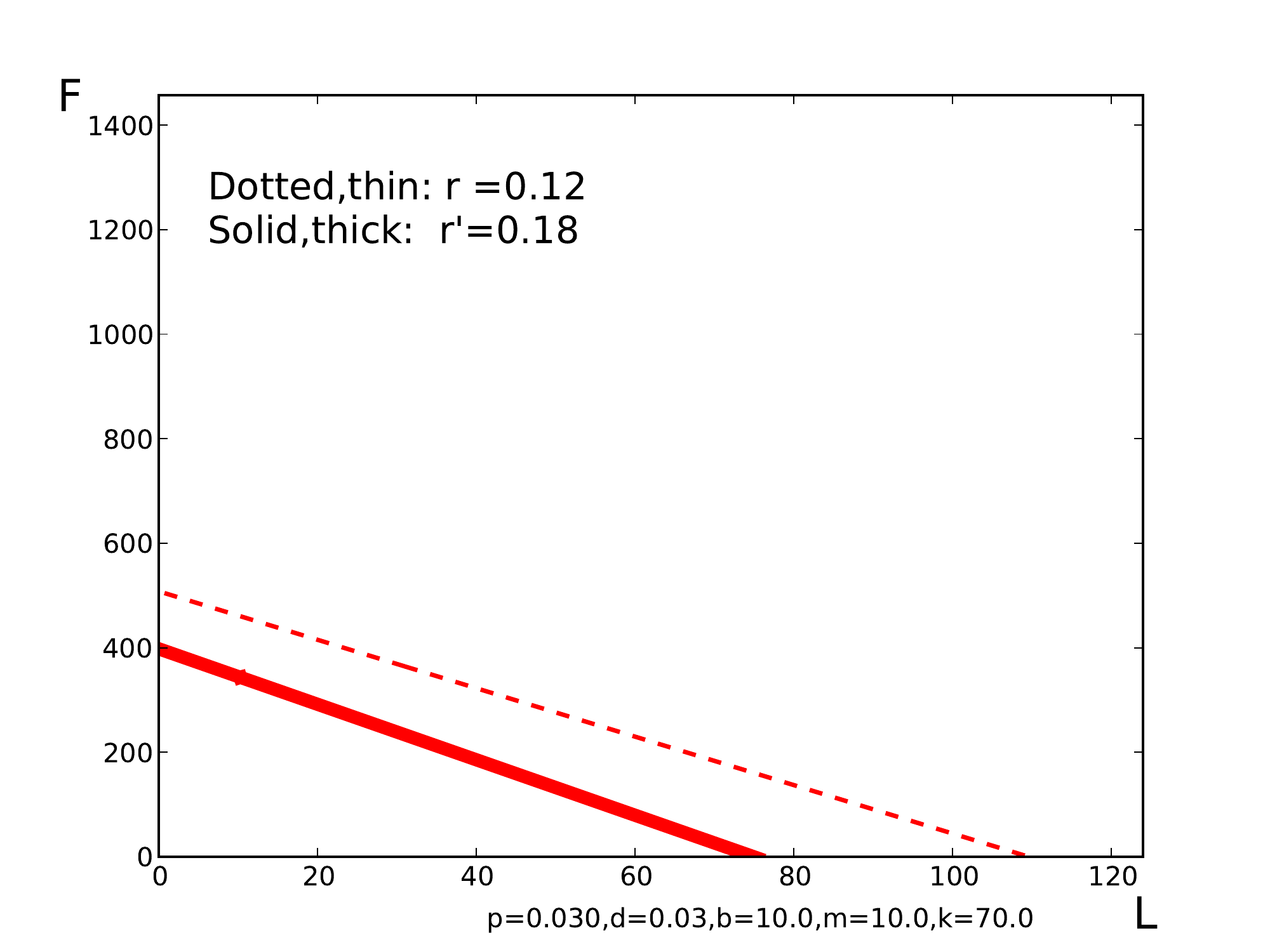}}\subfigure[]{\includegraphics[width=0.3\columnwidth]{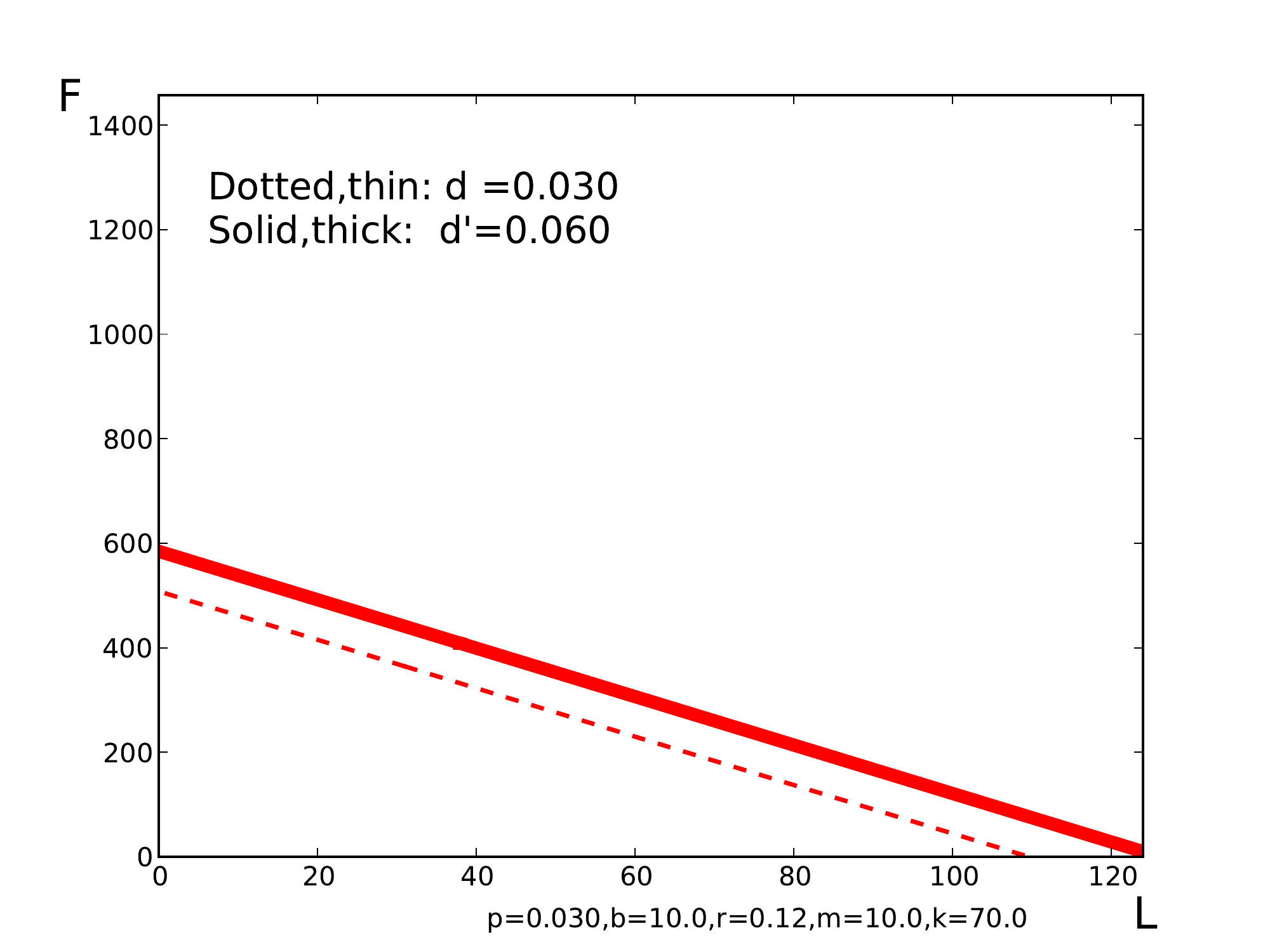}}
\par\end{centering}

\caption{\label{fig:The-effects-of-prd}The effects of the parameters $p$
(a), $r$ (b) and $d$ (c) on the dynamical system as seen through
the effect on the sink line. When $p$ or $r$ are increased the organizations
are able to grow faster, causing the sink line to move down, making
the existing CT measures no longer sufficient to neutralize some terrorist
organizations. In contrast, when $d$ is increased, the sink line
moves up because the organization is forced to replace more internal
loses to survive.}

\end{figure}

Suppose now that the model predicts that the given organization is
expected to grow further despite the current CT measures, and therefore
increased CT measures would be needed to defeat it. To see the effect
of additional CT measures, we need to examine how the dynamical system
changes in response to increases in the values of the parameters,
in particular, the parameters $b$ and $k$ which express the CT measures
directed at leaders and foot soldiers, respectively (Fig.~\ref{fig:The-effects-of-bk}).
It is also possible to affect the fate of the organization by influencing
the values of other parameters affecting its evolution, such as recruitment
and promotion (Fig.~\ref{fig:The-effects-of-prd}). In general, to
bring the terrorist organization under control it is necessary to
change the parameters individually or simultaneously so that the organization's
current state, $(L,F),$ is trapped under the sink line. An interesting
finding in this domain is that both $b$ and $k$ are equivalent in
the sense that both shift the sink link up in parallel (Fig.~\ref{fig:The-effects-of-bk}).

\section{Discussion}

\subsection{\label{sub:classification}Nascent terrorist organizations}

Recall that the sink line (Fig.~\ref{fig:Phase-Schematical}) distinguishes
two classes of terrorist organizations - those destined to be neutralized
and those that will continue growing indefinitely. Within the latter
group, another distinction is introduced by the trend line - a distinction
with significance to counter-terrorism efforts: organizations
lying to the left of it have different initial growth patterns compared
to those lying to the right (Fig.~\ref{fig:Phase-portrait}). The former
start with a large base of foot soldiers and a relatively small core
of leaders. In these organizations, $F$ may initially decline because
of CT, but the emergence of competent leaders would then start organizational
growth (e.g. Fig.~\ref{fig:Evolution1}(a)). In contrast, the latter
type of organizations start with a large pool of leaders but comparatively
few recruits. CT could decimate their leadership, but they would develop
a wide pool of foot soldiers, recover and grow (e.g. Fig.~\ref{fig:Evolution1}(b)).
Thus, all successful terrorist organizations may be classified as
either ``p-types'' (to the left of the trend line) or ``r-types''
(to the right of the trend line) in reference to the parameters $p$
of promotion and $r$ of recruitment. In p-type organizations early growth
occurs mainly through promotion of their foot soldiers to leaders, while in the
r-types mainly through recruitment of new foot soldiers.

This classification could be applied to many actual organizations.
For example, popular insurgencies are clearly p-type, while al-Qaeda's
history since the late 1990s closely follows the profile of an r-type:
Al-Qaeda may be said to have evolved through three stages: First,
a core of followers moved with bin Laden to Afghanistan. They were
well-trained but the organization had few followers in the wider world
\citep[for a history see][]{Wright06}. Then the attacks and counter-attacks
in the Fall of 2001 reduced the organization's presence in Afghanistan
leaving its operatives outside the country with few leaders or skills. Finally
the organization cultivated a wide international network of foot soldiers
but they were ill-trained as compared to their predecessors. This
description closely matches the profiles in Fig.~\ref{fig:Phase-portrait}
where r-type organizations start from a small well-trained core, move
toward a smaller ratio of leaders to foot soldiers and then grow through
recruitment. 

As was noted, nascent organizations tend towards the trend line, regardless
of how they started (Fig.~\ref{fig:Phase-portrait}). The slope of
this line is $\frac{r+\sqrt{r^{2}+4rmp}}{2p}$, and this number is
the long-term ratio between the number of foot soldiers and the number
of leaders. Notice that this formula implies that ratio is dependent
on just the parameters of growth - $r,m,p$ - and does not depend
on either $d$ or the CT measures $k,b$. This ratio is generally
not found in failing organizations, but is predicted to be ubiquitous
in successful organizations. It may be possible to estimate it by
capturing a division of an organization and it can help calculate
the model's parameters. However, it is important to note that $L$
includes not just commanding officers, but also any individuals with
substantially superior skills and experience. The existence of the
ratio is a prediction of the model, and if the other parameters are
known, it could be compared to empirical findings.

\subsection{\label{sub:Conditions-for-Victory}Conditions for Victory}

Recall, that the model indicates that all terrorist organizations
belong to one of three classes: r-types, p-types and organizations
that will be defeated. Each class exhibits characteristic changes
in its leaders, foot soldiers and strength ($L,F$ and $S$ resp.)
over time. This makes it possible to determine whether any given organization
belongs to the third class, i.e., to predict whether it would be defeated.

One finding is that if a terrorist organization weakens, i.e. shows
a decline in its strength $S$, it does not follow that it would be
defeated. Indeed, in some r-type organizations it is possible to observe
a temporary weakening of the organization and yet unless counter-terrorism
(CT) measures are increased, the organization would recover and grow
out of control (see Fig.~\ref{fig:Evolution1}(b)). Even a decline
in the leadership is not by itself sufficient to guarantee victory.
The underlying reason for this effect is out-of-control growth in
$F$, which would ultimately create a new generation of terrorist
leaders. Similarly, it is possible for an organization to experience
a decline in its pool of foot soldiers and yet recover. These cases
indicate that it is easy during a CT campaign to incorrectly perceive
apparent progress in reducing the organization as a sign of imminent
victory. 

Fortunately, under the model it is possible to identify reliable conditions
for victory over the organization (see the appendix for the proof):

\begin{enumerate}
\item \emph{For a p-type organization, it is impossible to have a decline
in strength $S$. If such a decline is made to happen, the organization
would be defeated.}
\item \emph{For an r-type organization, it is impossible to have a decline
in foot soldiers $F$. If such a decline is made to happen, the organization
would be defeated.}
\end{enumerate}
Consequently:

\begin{quote}
\emph{A terrorist organization would collapse if counter-terrorism
measures produce both: (1) a decline in its strength $S$ and (2)
a decline in its foot soldiers $F$.}
\end{quote}
In a notable contrast, declines in strength and the \emph{leaders}
are \emph{not} sufficient in all cases (see Fig.~\ref{fig:Evolution1}(b)).
To apply the theorem to an organization of an unknown type, one needs
merely to estimate whether the organization's pool of foot soldiers
and strength are declining\emph{}. The latter could be found indirectly
by looking at the quantity and quality of terrorist operations. It
is not necessary to know the model's parameters or changes in the
pool of leaders - the latter could even be increasing. Furthermore,
while it may take some time to determine whether $S$ and $F$ are
indeed declining, this time could be much shorter compared to the
lifetime of the organization. Therefore, the theorem suggests the
following two-step approach:

\begin{enumerate}
\item Estimate the scale of CT measures believed to be necessary to defeat
the organization. 
\item Measure the effect on $S$ and $F$. If they both declined, then sustain
the scale of operations (i.e. do not reduce $b$, $k$); Otherwise
an increase in CT measures would be necessary.
\end{enumerate}
The theorem and findings above give sufficient conditions for victory
but they do not characterize the only possible victory scenario. For
example, it is possible for an organization to see an increase in
its pool of foot soldiers $F$ yet ultimately collapse: these are
organizations that lie to the right of the trend line and just slightly
under the sink line. More generally, it should be remembered that
to prove the theorem it was necessary to use a simplified model of
a terrorist organization, as described in section \ref{sec:A-Mathematical-Model}.
Nevertheless, it is likely that some form of the theorem would remain
valid in complicated models because the model is built on fundamental
forces that are likely to be retained in these models.

\subsection{\label{sub:Stable-Equilibria}Stable Equilibria}

Recall that the model does not have a stable equilibrium (Fig.~\ref{fig:Phase-Schematical}).
Yet, in many practical cases, terrorist organization seem to reach
a stable equilibrium in terms of their structure and capabilities.
It is plausible that such stability is the result of a dynamic balance
between the growing terrorist organization and increasing CT measures
directed against it. Indeed, rather than staying constant numbers
like $b,k$, CT may actually grow when the organization presents more
of a threat%
\footnote{It would be a straightforward task to modify the model to incorporate
such a control-theoretic interaction, but the task is more properly
the subject of a follow-up study.%
}. Aside from CT, stability may be the result of organizations reaching
an external limit on their growth - a limit imposed by constraints
such as funding, training facilities or availability of recruits.
The case of funding could be modeled by assuming that the growth of
the organization slows as the organization approaches a maximum point,
($L_{max},F_{max}$).  Alternatively, it is quite possible and consistent
with the model that there would be a perception of stasis because
the organization is changing only slowly.

\section{\label{sec:Towards-an-Optimal}Counter-Terrorism Strategies}

Recall that the general counter-terrorism (CT) strategy in this model
is based on the location of the sink line, which we want to place above
the terrorist organization (in Fig.~\ref{fig:Phase-portrait}). To
implement this strategy, it is necessary first to calculate the model's
parameters for a given organization ($p,r,m,d$), and second, to determine
the efficacy of the current counter-terrorism measures ($b,k$). Then,
it remains {}``just'' to find the most efficient way of changing
those parameters so as to move the sink line into the desired location.
Let us now consider several strategic options.

\subsection{Targeting the leaders}

An important {}``counter terrorist dilemma'' \citep{Ganor05} is
whom to target primarily - the leaders or the foot soldiers. Foot
soldiers are an inviting target: not only do they do the vital grunt
work of terrorism, they also form the pool of potential leaders, and
thus their elimination does quiet but important damage to the future
of the organization. Moreover, in subsection \ref{sub:classification}
we saw that while an organization can recover from a decline in both
its strength and leadership pool, it cannot recover from declines
in both its strength and its foot soldiers pool. That finding does
not say that attacking leaders is unlikely to bring victory - indeed,
they form an important part of the organization's overall strength,
but it does suggest that a sustained campaign against an organization
is more likely to be successful when it includes an offensive against
its low-level personnel. Yet, it seems that the neutralization of
a terrorist leader would be more profitable since the leader is more
valuable to the organization than a foot soldier, and his or her loss
would inevitably result in command and control difficulties that may
even disrupt terrorist attacks.

When we use the model to address the problem quantitatively, we find
that the optimal strategy is actually dependent upon the organization,
that is to say the parameters $p,d,r,m$ (but not on $b,k$). For
example, for the parameter values used in the figures above, an increase
in $b$ gives a greater rise in the sink line than an equal increase
in $k$. Specifically, for those parameter values every two units
of $b$ are equivalent to about ten units of $k$. In general, when
$m,r,d$ are high but $p$ is low then attack on the leadership is
favored, while attack on the foot soldiers is best when $p$ is high
but $m,r,d$ are low - in agreement with intuition%
\footnote{Mathematically to obtain this result we first compute the derivatives
of the fixed point with respect to both $b$ and $k$, then project
them to the orthogonal to the sink link and then use an optimization
solver to find the parameter values of the model which maximize (and
minimize) the ratio between the lengths of the projections.%
}. In the first parameter range, foot soldiers are recruited so rapidly
that attacking them is futile, while in the second set leaders are
produced quickly so the only strategy is to fight the foot soldiers
to prevent them from becoming leaders. In any case, policy prescriptions
of this kind must be applied with consideration of counter-terrorism
capabilities and policy costs. Thus, while on paper a particular strategy
is better, the other strategy could be more feasible.

It is often argued that counter-terrorism policies have considerable
side effects. For instance, there is evidence that targeted assassinations
of leaders have led terrorist organizations to escalate, in what has
been called the ``boomerang effect'' \citep[p.125]{Crenshaw96}.
Fortunately, the model suggests that the policy maker has useful substitutes,
with possibly fewer policy side effects. As Fig.~\ref{fig:The-effects-of-prd}
shows, making recruitment ($r$) lower has an effect similar to increasing
$k$. Likewise, decreasing the rate of promotion to leadership ($p$)
can substitute for increasing $b$. This agrees with intuition: for
example, in the case of the foot soldiers, growth can be contained
either actively through e.g. arrests or proactively by slowing the
recruitment of new operatives (through e.g. attacks on recruitment
facilities or advocacy).

\subsection{\label{sub:Encouraging-desertion}Encouraging desertion}

Fatigue and attrition of personnel have been empirically found to
be an important effect in the evolution of terrorist organizations.
In interviews with captured or retired terrorists, they often complained
about the psychological stress of their past work, its moral contradictions,
and the isolation from relatives and friends \citep[Ch.6]{Horgan05}.
This is part of the reason why terrorist organizations cannot remain
inactive (as in a cease fire) for very long without experiencing irreplaceable
loss of personnel due to loss of motivation, and many organizations
even resort to coercion against desertion. Therefore, encouraging
operatives to leave through advocacy or amnesties may be an effective
counter-terrorism strategy. 

The model introduced here brings theoretical insight into this phenomenon.
One prediction of the model is that even if such desertion exceeds
recruitment (i.e. $d>r$) the organization would still sustain itself
as long as it has a sufficiently large rate of promotion ($p$) or
leaders of sufficiently high caliber ($m$). However, if $d$ is even
greater, namely, exceeds $d=\frac{1}{2}(r+r\sqrt{1+\frac{mp}{r}})$,
then the model predicts that the organization would be destroyed 
regardless of starting conditions, or counter-terrorism efforts ($b,k$). 

\begin{figure}[H]
\begin{centering}
\subfigure[]{\includegraphics[width=0.4\columnwidth]{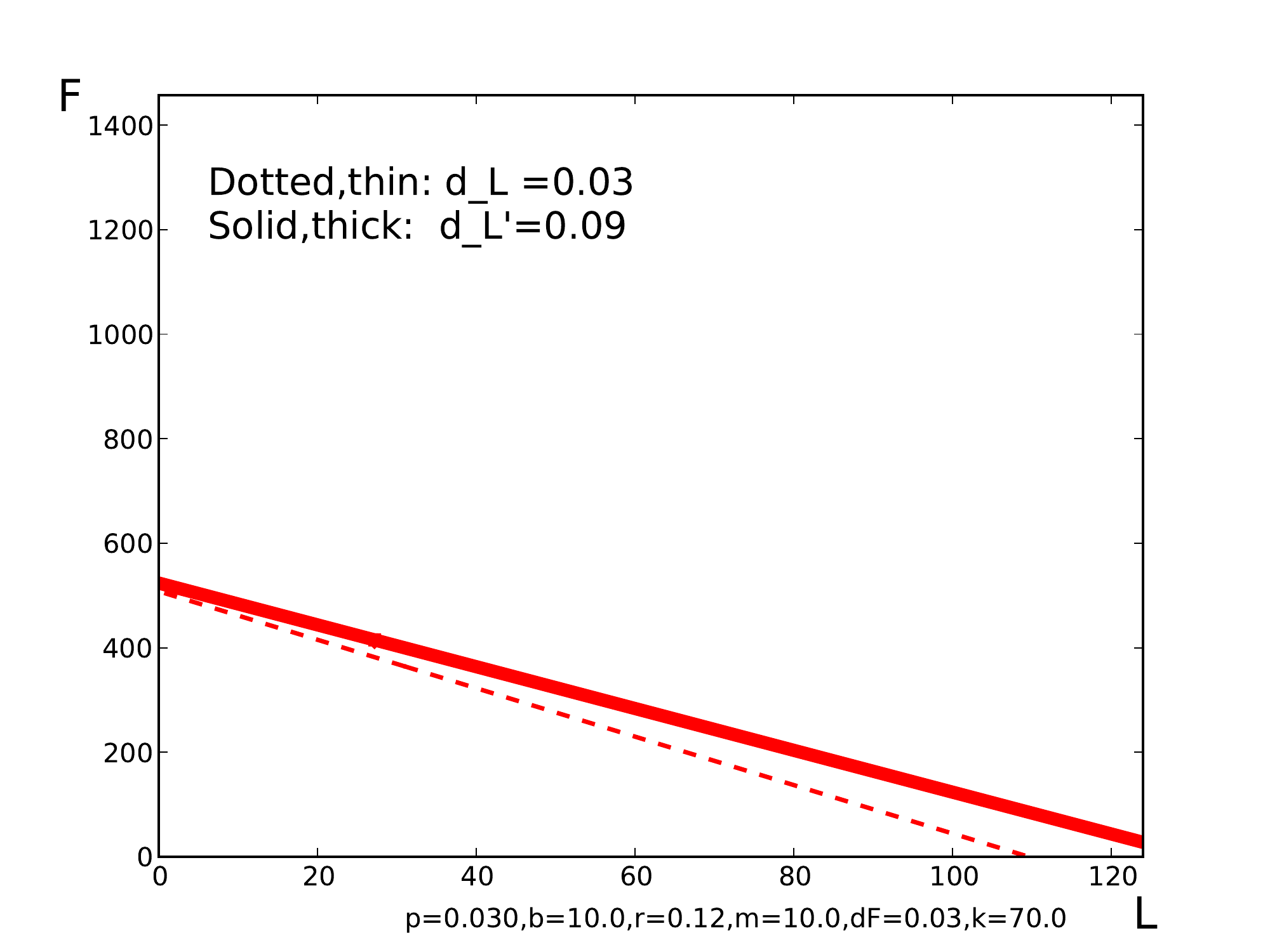}}\subfigure[]{\includegraphics[width=0.4\columnwidth]{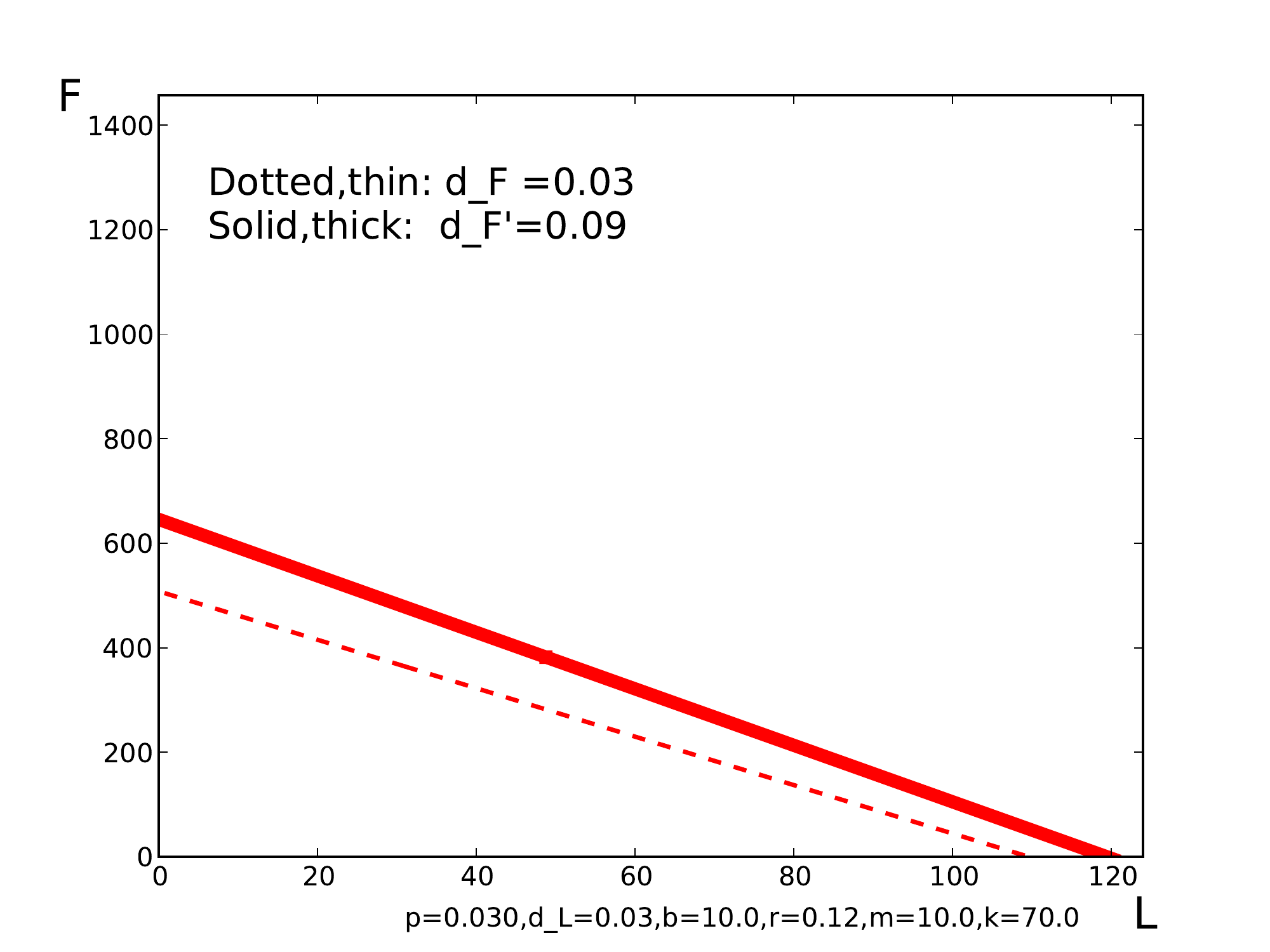}}
\caption{\label{fig:The-effects-of-d}The effects of $d_{L}$ (a) and $d_{F}$ (b) 
on the dynamical system, as seen through the effect on the sink
line. As the desertion rates increase, the sink line moves up and
its slope changes, thus trapping additional terrorist organizations. }

\par\end{centering}
\end{figure}

Organizations with lower $d$ are, of course, also effected by desertion.
Earlier, in Fig.~\ref{fig:The-effects-of-prd} we saw how increasing
$d$ raises up the sink line. To see the phenomenon in more detail,
we replaced $d$ by two (not necessarily equal) parameters $d_{L}$
and $d_{F}$ for the desertion of $L$ and $F$, respectively. The
two parameters change the slope of the sink line: increasing $d_{L}$
flattens it, while increasing $d_{F}$ makes it more steep (Fig.~\ref{fig:The-effects-of-d}).
Therefore, increasing $d_{L}$ could be a particularly effective strategy
against nascent r-type organizations, while increasing $d_{F}$ could
be effective against the nascent p-types.

\subsection{\label{sub:Minimization-of-Strength}Minimization of Strength $S$}

Counter terrorism (CT) is often the problem of resource allocation
among competing strategies. Therefore, suppose that resources have
become available towards a CT operation against the terrorist organization.
Namely, suppose we can remove leaders and operatives in a single blow
(unlike the parameters $b,k$ in the model which take a gradual toll).
A reasonable approach to allocating these resources efficiently would
be to divide them between operations targeting the leadership and
those targeting the foot soldiers in such a way that the terrorist
organization's strength $S$ is minimized%
\footnote{Mathematically, this would be two variable minimization of $S$ constrained
by a budget.%
}. However, by some simplified economic analysis, it is possible to
show (see appendix) that this counter-terrorism strategy is in general \emph{suboptimal}.
Instead, for a truly effective resource allocation, it is necessary
to consider the dynamics of the organization being targeted and the
true optimum may be considerably different. For example, when the
ratio of promotion to recruitment is relatively large (i.e. $\frac{p}{r}\gg0$),
then the optimum shifts increasingly towards attacking the foot soldiers
since they become much harder to replace than leaders.

On an intuitive level, the reason why the strategy is suboptimal is
because often, the losses we can inflict most effectively on the organization
are precisely those losses that the organization can restore most
easily. Hence, in the long-term a strategy targeting strength $S$
would be ineffective. Instead, when making a CT strategy it would
be valuable to understand the target organization's dynamics, and
in particular, to build a dynamical model. Such a model would help
because it can identify an organization's unique set of vulnerabilities
due to the composition of its human capital and its properties as
a dynamical system.

\section{Conclusions}

Much of the benefit of mathematical models is due to their ability
to elucidate the logical implications of empirical knowledge that
was used to construct the model. Thus, whereas the empirical facts
used to construct the models should be uncontroversial, their conclusions
should offer new insights. The model proposed here is a very simplified
description of real terrorist organizations. Despite its simplicity,
it leads to many plausible predictions and policy recommendations.
Indeed, the simplicity of the model is crucial to making the model
useful. More detailed models of this kind could provide unparalleled
insights into counter-terrorism policies and the dynamics of terrorism.
\\

\noindent\textbf{Acknowledgments}

I would like to thank Steven Strogatz, Michael Genkin, and Peter Katzenstein
for help and commentary, Richard Durrett for unconditional support,
R.A. for infinite patience and my colleagues at the Center for Applied Mathematics for encouragement.
\\

\noindent\textbf{Note}
An earlier shorter version of this paper appeared in \textit{Studies in Conflict and Terrorism},
volume $32$, number $1$, $2009$.

\appendix

\section{Appendix}

The original differential equations are:

\begin{eqnarray}
\frac{dL}{dt} & = & pF-d\, L-b\label{eq:L}\\
\frac{dF}{dt} & = & r(\underset{S}{\underbrace{mL+F}})-d\, F-k\label{eq:F}\end{eqnarray}
If we wished to incorporate the drain of the foot soldiers due to
promotion ($-pF$) in Eqn.\ref{eq:F}, then we could adjust the original
parameters by the transformation $r\to r-p$ and $m\to\frac{rm}{r-p}$.
However, this would affect some of the analysis below, because for
$r<p$ it would not longer be the case that $r>0$, even though $rm>0$
would still hold true. Alternatively, we could change the internal
losses parameter for foot soldiers : $d_{F}\to d_{F}+p$ and break
the condition $d_{F}=d_{L}$.

The linearity of the system of differential equations makes it possible
to analyze the solutions in great detail by purely analytic means.
The fixed point is at:\begin{eqnarray}
L_{*} & = & \frac{kp-b(r-d)}{d(r-d)+rmp}\quad F_{*}=\frac{kd+rmb}{d(r-d)+rmp}\label{eq:fixedPt}\end{eqnarray}
The eigenvalues at the fixed point are\begin{eqnarray}
\lambda_{1,2} & = & \frac{r-2d\pm\sqrt{(r-2d)^{2}+4\left(rmp+d(r-d)\right)}}{2}\label{eq:fixedPtEvals}\end{eqnarray}
From Eqn.(\ref{eq:fixedPtEvals}), the fixed point is a saddle when
$rmp+d(r-d)>0$, i.e. $\frac{r-\sqrt{r^{2}+4rmp}}{2}<d<\frac{r+\sqrt{r^{2}+4rmp}}{2}$
(physically, the lower bound on $d$ is $0$). The saddle becomes
a sink if $r<2d$ and $rmp+d(r-d)<0$. By Eqn.(\ref{eq:fixedPt}),
this automatically gives $F_{*}<0$, i.e. the organization is destroyed%
\footnote{Of course, the dynamical system is unrealistic once either $F$ or
$L$ fall through zero. However, by the logic of the model, once $F$
reaches zero, the organization is doomed because it lacks a pool of
foot soldiers from which to rebuild inevitable losses in its leaders.%
}. It is impossible to obtain either a source because it requires $r-2d>0$
and $rmp+d(r-d)<0$, but the latter implies $d>r,$ and so $r-2d>0$
is impossible; or any type of spiral because $(r-2d)^{2}+4(rmp+d(r-d))<0$
is algebraically impossible%
\footnote{The degenerate case of $\lambda=0$ has probability zero, and is not
discussed.%
}. It is also interesting to find the eigenvectors because they give
the directions of the sink and trend lines:\begin{equation}
e_{1,2}=\left(\begin{array}{c}
2p\\
r\pm\sqrt{r^{2}+4rmp}\end{array}\right)\label{eq:slopes}\end{equation}
We see that the slope of $e_{2}$, which is also the slope of the
sink line - the stable manifold - is negative. Therefore, we conclude
that the stable manifold encloses, together with the axes, the region
of neutralized organizations. Concurrently, the slope of $e_{1}$
- the trend line i.e. the unstable manifold - is positive. Thus, the
top half of the stable separatrix would point away from the axes,
and gives the growth trend of all non-neutralized organizations ($\frac{\Delta F}{\Delta L}=\frac{r+\sqrt{r^{2}+4rmp}}{2p}$).

\subsection{Proof of the theorem}

Recall, we wish to show that a terrorist organization that experiences
both a decline in its strength and a decline in the number of its
foot soldiers will be destroyed. The proof rests on two claims: First,
a p-type organization cannot experience a decline in strength, and
second, an r-type organization cannot experience a decrease in $F$
(for a graphic illustration see Fig.~\ref{fig:Isostrength-Proof} below).
Thus, both a decline in strength and a decline in the number of foot
soldiers cannot both occur in an r-type organization nor can they both
occur in a p-type organization. Hence, such a situation can only occur
in the region of defeated organizations.

\begin{wrapfigure}{r}{0.6\textwidth}
\begin{centering}
\includegraphics[width=0.6\textwidth]{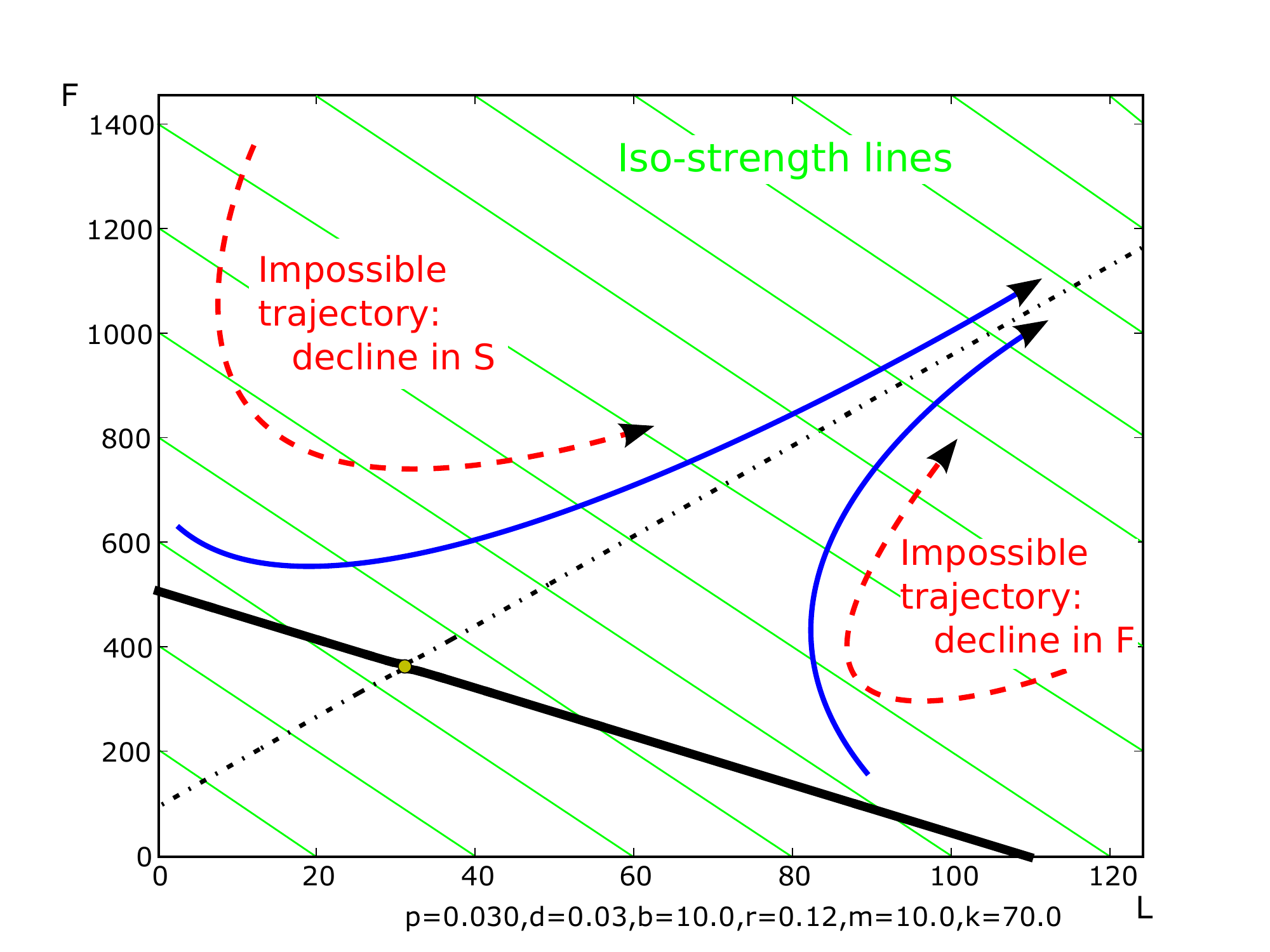}
\caption{\label{fig:Isostrength-Proof}The phase plane with possible (solid blue)
and impossible (dashed red) trajectories, and lines of equal organization
strength (green). Because orbits of the p-type must experience an
increase in strength $S$, the left red line cannot be an orbit. Also,
r-type orbits must experience an increase in $F$, and so the right
red line cannot be an orbit either.}
\par\end{centering}
\end{wrapfigure}

As to the first claim, we begin by showing that the slope of the sink
line is always greater than the slope of the iso-strength lines ($=-m$).
By Eqn. (\ref{eq:slopes}) the slope is $\frac{r-\sqrt{r^{2}+4rmp}}{2p}=-m\frac{2}{1+\sqrt{1+4\frac{mp}{r}}}>-m$.
Therefore, the flow \emph{down} the sink line has $\frac{dS}{dt}>0$
(Down is the left-to-right flow in the figure). Now, we will show
that in a p-type organization, the flow must experience an even greater
increase in strength. Let $A$ be the matrix of the dynamical system
about the equilibrium point and let the state of the terrorist organization
be $(L,F)=d_{1}e_{1}+d_{2}e_{2}$ where $e_{1},e_{2}$ are the distinct
eigenvectors corresponding to the eigenvalues $\lambda_{1},\lambda_{2}$.
Consideration of the directions of the vectors (Eqn.(\ref{eq:slopes}))
shows that for a p-type organization, $d_{1}>0$ and $d_{2}<0$. The
direction of flow is therefore $d_{1}\lambda_{1}e_{1}+d_{2}\lambda_{2}e_{2}$.
Notice that $\lambda_{1}>0,$$\lambda_{2}<0$, and so the flow has
a positive component $(=d_{1}\lambda_{1})$ in the $e_{1}$ direction
(i.e. up the trend line). Since the flow along $e_{1}$ experiences
an increase in both $L$ and $F$, it must experience an increase
in strength. Consequently, a p-type organization must have $\frac{dS}{dt}$
which is even more positive than the flow along the sink line (where
$d_{1}=0$). Thus, $\frac{dS}{dt}>0$ for p-types. 

As to the second claim, note that r-type organizations have $d_{1}>0$
and $d_{2}>0$. Moreover, in an r-type organization, the flow $d_{1}\lambda_{1}e_{1}+d_{2}\lambda_{2}e_{2}$
has $\frac{dF}{dt}$ greater than for the flow \emph{up} the right
side of the sink line (right-to-left in the figure): the reason is
that $e_{1}$ points in the direction of increasing $F$ and while
in an r-type $d_{1}>0$, along the sink line $d_{1}=0$. The flow
up the sink line has $\frac{dF}{dt}>0$, and so $\frac{dF}{dt}>0$
in an r-type organization. In sum, $\frac{dS}{dt}<0$ simultaneously
with $\frac{dF}{dt}<0$ can only occur in the region $d_{1}<0$ -
the region of defeated organizations. QED.

\subsection{Concrete Example of Strength Minimization}

\begin{figure}[H]
\begin{centering}
\subfigure[]{\includegraphics[width=0.4\columnwidth]{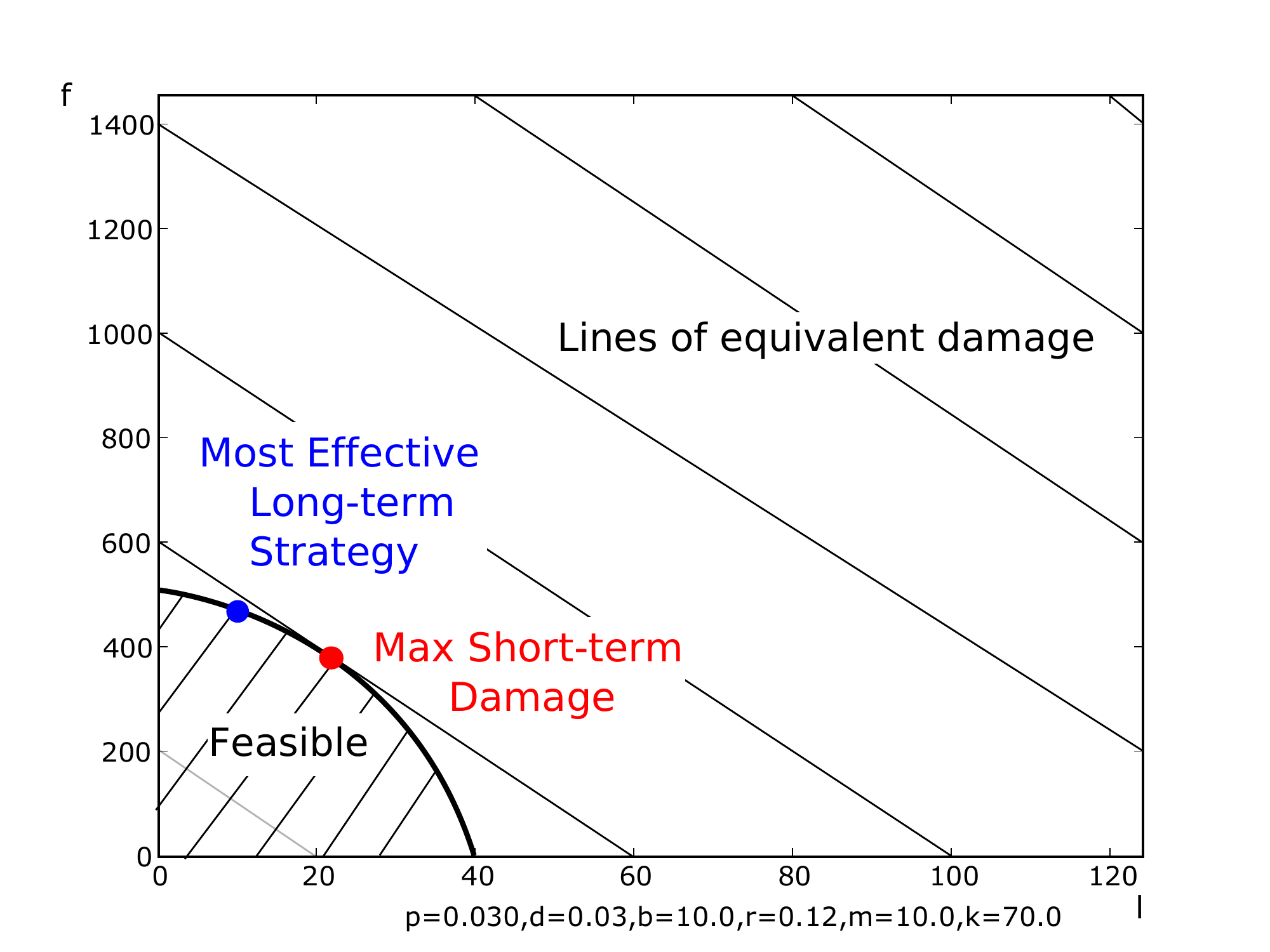}}\subfigure[]{\includegraphics[width=0.4\columnwidth]{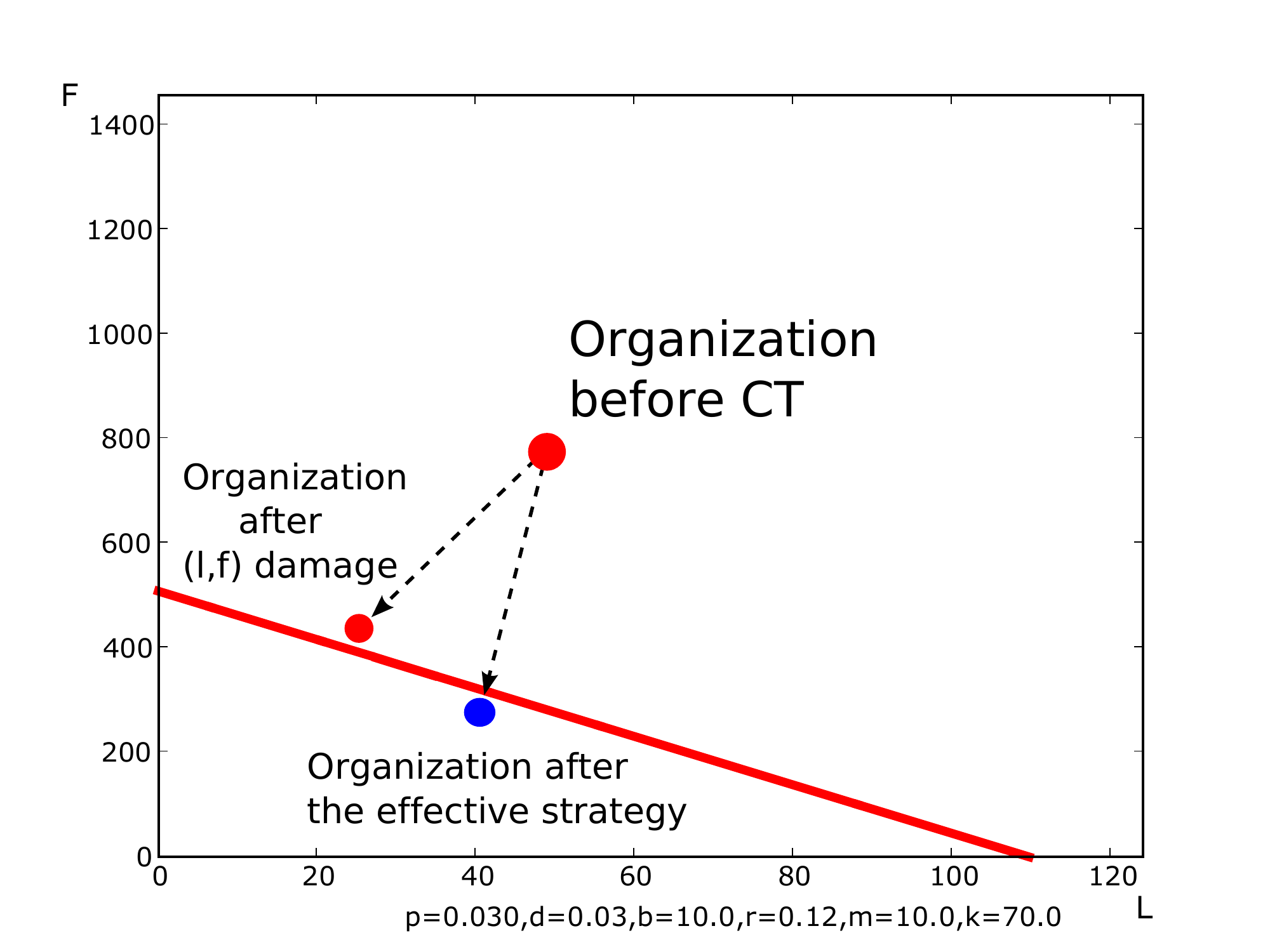}}
\par\end{centering}
\caption{\label{fig:CT-optimum}Graphical calculation of optimal budget allocation
(a) and contrast between minimization of $S$ and the actual optimum
(b). In (a), the optimal choice of ($l,m$) is given by the point
of tangency between the feasible region and the lines of constant
$S$. In (b), the red line is the sink line. The minimization of $S$
through the removal of about $20$ leaders and $400$ foot soldiers
would not bring the organization under the sink line, but a different
(still feasible) strategy would.}
\end{figure}

In subsection \ref{sub:Minimization-of-Strength} we claimed that
the task of minimizing $S$ is different from the optimal counter-terrorism
strategy. Here is a concrete example that quantitatively illustrates
this point. Suppose a resource budget $B$ is to be allocated between
fighting the leadership and fighting the foot soldiers, and furthermore,
that the cost of removing $l$ leaders and $f$ foot soldiers, respectively,
is a typical convex function: $c_{1}l^{\sigma}+c_{2}f^{\sigma}$ ($c_{1}$
and $c_{2}$ are some positive constants and $\sigma>1$)%
\footnote{$\sigma>1$ because e.g. once the first say $20$ easy targets are
neutralized, it becomes harder to find and neutralize the next $20$
(the law of diminishing returns.) In any case the discussion makes clear
that for most cost functions the suggested optimum would be different 
from the true optimum.%
}. Notice that whereas uppercase letters $L,F$ indicate the number
of leaders and foot soldiers, respectively, we use lowercase $l,f$
to indicate the number to be \emph{removed}. The optimal values of
$l$ and $f$ can be easily found graphically using the standard procedure
in constrained optimization: the optimum is the point of tangency
between the curve $B=c_{1}l^{\sigma}+c_{2}f^{\sigma}$ and the lines
of constant difference in strength: $\Delta S=ml+f=constant$ (Fig.\ref{fig:CT-optimum}(a)).
However, as illustrated in Fig.~\ref{fig:CT-optimum}(b), if such a
strategy is followed, the terrorist organization may still remain
out of control. It is preferable to choose a different strategy -
in the example it is the strategy that focuses more on attacking the foot
soldiers and thus brings the organization under the sink line (red
line), even though the $\Delta S$ is not as large. In general, the
difference between the strategies is represented by the difference
between the slope of the sink line and the slopes of the lines of
equivalent damage to strength. The latter always have slope $-m$
while the former becomes arbitrarily flat as $\frac{p}{r}\to\infty$.

\bibliographystyle{chicagoa}
\bibliography{terror,math}

\end{document}